\documentclass[11pt]{amsart}

\usepackage{amssymb}

\usepackage{graphicx}

\usepackage{amsmath}

\usepackage{verbatim}

\makeindex

\input xy
\xyoption{all}

\newcommand{\prskip}{\vspace{8pt}} 
\newcommand{\thmskip}{\vspace{10pt}} 
\newcommand{\sectskip}{\vspace{50pt}} 
\newcommand{\introskip}{\vspace{25pt}} 

\begin{document}

\date{15 June, 2020}
\title{A computing strategy and programs to \\resolve the Gerstenhaber Problem for\\ commuting triples of matrices}
\vspace{3mm}

\author{John Holbrook {\Small{and}} Kevin C.\ O'Meara}

\maketitle
\begin{abstract} We describe a {\Small{MATLAB}} program that could produce a negative answer to the Gerstenhaber Problem by the construction of three commuting $n \times n$ matrices $A,B,C$ over a field $F$ such that the subalgebra $F[A,B,C]$ they generate has dimension greater than $n$. This problem has remained open for nearly 60 years, following Gerstenhaber's surprising result (Annals Math.) that $\dim F[A,B] \le n$ for any two commuting matrices $A,B$. The property fails for four or more commuting matrices. We also make the {\Small{MATLAB}} files freely available.
\end{abstract}
\vspace{5mm}

Easily-stated, longstanding problems have always fascinated mathematicians (think of Fermat's Last Theorem).  The Gerstenhaber Problem (GP) for commuting matrices fits into this category. But to our knowledge, our project is the first serious attempt to use a computer to provide a negative answer, despite very strong evidence that this is the case (see, for example, \cite{HO}, \cite{O}). Recently, the second author showed in \cite{O} that the GP is Turing computable, marking a huge reduction in the complexity of the problem. Moreover, he showed it suffices to answer the GP in the case of the fields $F = \mathbb{Z}/p$ of integers mod $p$ for primes $p$. Thus the GP is ideally suited to a computer search. We outline a computing strategy that the authors have developed during the last few years. We will also make the {\Small{MATLAB}} files freely available upon request to the authors using the e-mail
\begin{center}
                \texttt{GPmfiles@gmail.com}
\end{center}
so that others can join in this safari. Graduate students with a knowledge of basic algebra and {\Small{MATLAB}} could find this an attractive project (and possible path to early fame!). Particularly those with access to fast computers and parallel computing.
\prskip

The authors are not experts in computing. One important goal of our paper is to encourage some of the very talented young folk out there, those who really are experts in computing, to join in. Undoubtedly fresh ideas would then flow. We chose to write our programs in the programming language {\Small{MATLAB}}, which was introduced in the early 1980's with particular emphasis on efficient handling of operations with large arrays (matrices being a special case). Nowadays  {\Small{MATLAB}} does so much more, of course, but it is still generally regarded as the best programming language for numerical linear algebra calculations with very large matrices. So many of our subroutines rely heavily on the basic ``calculus'' of matrices --- row operations and echelon forms (in our case, essentially for integer matrices, but at times of the order of $10,000 \times 20,000$). Mod $p$ arithmetic does not get around this. Some readers may prefer using a different programming language.
\prskip

The GP makes sense over any field $F$. However, the difficulties of quickly producing commuting triples and accurately calculating dimension depend on $F$. Even with the rational field $\mathbb{Q}$, where the GP can be reduced to integer matrices after clearing fractions, these calculations are no picnic. It is inherent in the GP that large powers of matrices (maybe the $50^{th}$ power of a $100 \times 100$ matrix) need to be considered. But then we quickly face roundoff errors. This will cast serious doubts on the accuracy of dimension (which is a highly discontinuous function of the matrix entries). Within {\Small{MATLAB}}, where the standard precision is 16 places, one can opt for higher precision by using variable precision arithmetic, vpa, to any specified degree of accuracy. But this dramatically slows calculations (and we expect to have to run hundreds of millions of trials to flush out a a counterexample for the GP).\footnote{\ This frustration with arithmetic blowouts was a strong motivation to go to mod $p$ calculations, and in turn a factor in leading the second author to the Turing computability of the GP. So in a sense, here a computer formulation led to an important conjecture and later proof.} On the other hand, matrix calculations done mod $p$ (even for a 4-digit prime) are fast and accurate.
\prskip

A crucial tool in our approach is the use of the Weyr canonical form for matrices, which is far better suited to commuting problems than its cousin, the Jordan form. In fact, it would be extremely difficult to formulate our approach in terms of the Jordan form (because of a simultaneous triangularization result mentioned in Section 2). In Section 2, we record some of the properties of the Weyr form, which are then built into our computer programs. Essentially these allow a recursive approach, starting with a pair of smaller-sized commuting matrices, to both the rapid construction of random commuting triples of large matrices and calculation of the dimension of the subalgebra they generate.
\prskip

It appears that the choice of settings for parameters in our programs, particularly the sparseness of the matrix entries, will be an important factor in producing some $\dim F[A,B,C] > n$, \,a case which the authors informally term a ``\emph{Eureka}''. We know Eurekas exist with commuting quadruples of matrices, but surprisingly they don't occur that frequently in computer calculations without carefully setting the parameters. We use such productive settings in the commuting quadruples case (for the same matrix size $n$) as a guide to how best to set them for commuting triples.
\prskip

Most of our testing to date has been done on just laptops, and working with commuting triples of matrices of size ranging from $15 \times 15$ to $50 \times 50$. The three sample trials given near the end of the paper indicate the sort of thing we have done. (Our suggestion to a potential user of our programs, as a first step to understanding our approach, is to re-run these sample trials.) However, the authors feel that the Loch Ness monster probably lives in deeper water, closer to $100 \times 100$. In that range, even more care is needed with sparsity settings because the number of fully completed trials will be small without very sparse matrices.
\prskip

There is a broader reason for publicising our approach. Many of the techniques can be adapted to other problems in commutative algebras of matrices, one demonstration of which was given in the 2017 paper \cite{OW} by O'Meara and Watanabe.  Although the programs in the present paper use mod $p$ arithmetic, it would be a simple matter for  readers to change the arithmetic to suit their needs (e.g.\ real or complex arithmetic). With readers in mind who might wish to modify our computer programs, perhaps using the languages {\Small{MAGMA}}  or {\Small{GAP}}, we have made the description of the programs given here (and documentation within the code) a high priority.
\introskip

\section{A brief history of the Gerstenhaber Problem}
\thmskip

The following is just a bare sketch of some of the history behind the GP. Much fuller accounts can be found in the references at the end of our paper. As for the authors' own accounts, Chapter 5 of \cite{ATLA} is devoted to Gerstenhaber's original theorem while Chapter 7 develops from scratch all the algebraic geometry necessary to understand how the powerful techniques of that area (irreducible varieties, algebraic geometry dimension etc) have impacted the GP. The paper \cite{HO} devotes many pages to the history and recent developments.
\prskip

In 1961 Gerstenhaber \cite{G1} established his (most surprising) theorem that \linebreak $\dim F[A,B] \le n$ for all commuting $n \times n$  matrices $A,B$ over any field $F$. (For a single ``commuting'' matrix $A$ this is an immediate consequence of the Cayley-Hamilton theorem.) His proof used algebraic geometry, but later purely linear-algebraic proofs were also found by Barr\'{\i}a-Halmos \cite{B1} and Laffey-Lazarus \cite{LL}. The earlier Holbrook paper \cite{H} may be mentioned in this connection, but it deals only with complex matrices. Whether Gerstenhaber's result extends to three commuting matrices $A,B,C$ is what is known as the \textbf{\emph{Gerstenhaber Problem}} \textbf{(GP)}. That is, must $\dim F[A,B,C] \le n$ for all $n \times n$ commuting matrices $A,B,C$? Recall $F[A,B,C]$ is the smallest(unital) subalgebra of the matrix algebra $M_n(F)$ of all $n \times n$ matrices over $F$ that contains $A,B,C$. Thus here, $F[A,B,C]$ is the space of all commuting polynomials in $A,B,C$ (including the scalar matrices $cI$). Without commutativity we can easily have $F[A,B] = M_n(F)$, giving $\dim F[A,B] = n^2$, the maximum possible.
\prskip

Historically, two of the most significant papers related to Gerstenhaber's theorem are the 1955 Motzkin-Taussky paper \cite{M1} and the 1992 Robert Guralnick paper \cite{G2}. The former established that the variety $\mathcal{C}(2,n)$ of commuting pairs of $n \times n$ matrices is always irreducible. Note this paper predates Gerstenhaber's paper \cite {G1} but Motzkin and Taussky in \cite{M1} seemed unaware of Gerstenhaber's result to come later. (Also Gerstenhaber doesn't mention \cite{M1}, although in his paper he establishes the Motzkin-Taussky irreducibility result.)  Guralnick in his paper simplified the Motzkin-Taussky argument and showed how quickly Gerstenhaber's theorem follows from the Motzkin-Taussky theorem.\footnote{\ This in no way lessens the credit due to Gerstenhaber for having the mathematical insight that the dimension bound for commuting pairs might be the same as for a single matrix --- the result is so unexpected and few would have even conjectured it.}  Moreover, he showed that the variety of commuting triples of matrices need not be irreducible, and therefore that method of attack could not answer the GP. It also gave weight to the possibility of a negative answer to the GP.
\prskip

The current state-of-play is there are no examples known where $\dim F[A,B,C] > n$. However, it is known that the GP has a positive answer for $n = 1,2, \ldots ,11$ when $F$ is of characteristic 0. The proofs are long and intricate, and use a mix of algebraic geometry and matrix canonical forms from linear algebra (see for instance Klemen {\v{S}}ivic's papers).
\prskip

By an algebraic geometry argument, it is known that, for a fixed matrix size $n$, $\dim F[A,B,C] \le n$ over an algebraically closed field $F$ (and therefore for subfields thereof) whenever the variety $\mathcal{C}(3,n)$ of all commuting triples of $n \times n$ matrices is irreducible. This is true in characteristic 0 for $n < 12$ (Klemen {\v{S}}ivic  and others) but not when $n > 28$ \,(Robert Guralnick \cite{G2}, Holbrook-Omladi\v{c} \cite{H3}). What happens in the integer interval [12,\,28] remains open.
\prskip

The linear-algebraic proofs of Gerstenhaber's theorem proceed inductively on the matrix size $n$, and centre on a fixed form spanning set for $F[A,B]$ of $n$ words in $A,B$ depending only on the Jordan (or Weyr) structure of the first matrix $A$. However, such fixed spanning sets for commuting triples do not exist. See \cite{O} for a discussion of this.
\prskip

Based on the above history (particularly the failure of the Motzkin-Taussky theorem to extend to commuting triples of matrices, and likewise the linear-algebraic techniques), and their own experiments, the authors firmly believe the GP will turn out to have a negative answer. The question is whether a specific counter-example is accessible with current computing techniques.
\thmskip

\section{The central role played by the Weyr form}
\thmskip

By a standard argument in linear algebra (using the generalized eigenspace decomposition), the GP reduces to considering commuting triples $A,B,C$ of \textbf{nilpotent} matrices, that is, some power of each matrix is the zero matrix. The smallest such power is called the \textbf{nilpotent index} of the matrix (note by the Cayley-Hamilton theorem, the nilpotent index of an $n \times n$ matrix $X \in M_n(F)$ can't exceed $n$, although always $X^n = 0$.) When $F$ is algebraically closed, the nilpotent matrices are exactly those whose only eigenvalue is 0. It is also well-known that, over an algebraically closed field, commuting matrices $X_1,X_2,\ldots, X_k$ can be simultaneously triangularised, that is, for some invertible matrix $P$, all the $P^{-1}X_iP$ are upper triangular (in fact must then be strictly upper triangular if the matrices are nilpotent). But as was shown by O'Meara and Vinsonhaler \cite{OV} in 2006 (see also Theorem 2.3.5 in \cite{ATLA} for an account), we can do much better than this:
\thmskip

\noindent {\Small{\textbf{CRITICAL PROPERTY FOR OUR APPROACH.}}}\, \emph{If $A_1,A_2,\ldots,A_k$ are commuting $n \times n$ matrices over an algebraically closed field $F$, then there is a similarity transformation that puts $A_1$ in Weyr form and simultaneously puts $A_2,\ldots,A_k$ in upper triangular form. (The algebraically closed assumption can be dropped if the matrices involved have all their eigenvalues in F, such as with nilpotent matrices.)}
\prskip

\noindent This fact applied to three commuting matrices $A,B,C$ is at the heart of our computing algorithm. It fails for the Jordan form.  Of course, it is enough for the purposes of the GP to establish whether $\dim F[A,B,C] \le n$ after such a simultaneous similarity transformation, because the transformed algebra is isomorphic to the original. In summary, and henceforth assumed:
\prskip

\noindent \textbf{The Reduction}.  {\emph{The commuting triple $A,B,C$ takes the form $W,K,L$ where $W$ is a nilpotent Weyr matrix and $K,L$ are strictly upper triangular}}.
\prskip

So what is the Weyr form? The form was introduced by the Czech mathematician Eduard Weyr in 1885, but then largely forgotten. It has been periodically rediscovered over the years by a handful of people (including O'Meara and Vinsonhaler), and now appears finally to be gaining traction. The 2nd edition \cite{HJ} of the Horn and Johnson ``\emph{Matrix Analysis}'' (2013)  covers the Weyr form in Chapter 3. The monograph \cite{ATLA} by O'Meara, Clark, and Vinsonhaler (2011) gives a comprehensive account of the form.  We will restrict our brief discussion of the Weyr form to nilpotent matrices.
\prskip

The best view of  a nilpotent Weyr matrix $W$ is that it is a blocked matrix form of a basic $r \times r$ nilpotent Jordan matrix
\[
    \renewcommand{\arraystretch}{1.2}
    J \ \ = \ \       \left[ \begin{array}{ccccccc}
           0  &  1       &       &       &        &     &     \\
              & 0        &   1   &       &        &     &     \\
              &          &\ddots  &       &        &     &     \\
              &          &       &       &        &   0 &   1 \\
              &          &       &       &        &     &   0
                 \end{array}   \right].
\]
In the simplest case where $n = dr$ for some specified positive integers $d,r$, the Weyr analogue is the $r \times r$ blocked matrix with $d \times d$ blocks
\[
    \renewcommand{\arraystretch}{1.2}
    W \ \ = \ \       \left[ \begin{array}{ccccccc}
           0  &  I       &       &       &        &     &     \\
              & 0        &   I   &       &        &     &     \\
              &          &\ddots  &       &        &     &     \\
              &          &       &       &        &   0 &   I \\
              &          &       &       &        &     &   0
                 \end{array}   \right].
\]
Here $0$ and $I$ are the $d \times d$ zero matrix and identity matrix respectively. In general, the blocking can be in terms of a specified partition $n = n_1 + n_2 + \cdots + n_r$ of $n$, with $n_1 \ge n_2 \ge \cdots \ge n_r$, that determines the sizes of the diagonal blocks. We then call $(n_1,n_2,\ldots,n_r)$ the  \textbf{Weyr structure} of $W$. The structure is said to be \textbf{homogeneous} if $n_1 = n_2 = \cdots = n_r$. The first superdiagonal blocks in the general case take the form
\[
      W_{j,j+1} \ = \    \left[\begin{array}{c}
                            I_j \\
                            0
                            \end{array}\right]
\]
where $I_j$ is the $n_{j+1} \times n_{j+1}$ identity matrix and $0$ is the $(n_j - n_{j+1}) \times n_{j+1}$ zero matrix (equivalently, $W_{j,j+1}$ is a reduced row-echelon matrix of full column rank). All other blocks are zero.  A picture is worth a thousand words:
\[
\renewcommand{\arraystretch}{1.1}
            W \ = \  \left[ \begin{array}{r r r |r r |r }
    0 & 0  & 0  &  1  & 0  & 0    \\
    0 & 0  & 0  &  0  & 1  & 0    \\
    0 & 0  & 0  &  0  & 0  & 0    \\  \cline{1-6}
      &    &    &  0  & 0  & 1    \\
      &    &    &  0  & 0  & 0    \\  \cline{4-6}
  \multicolumn{5}{c |}{ } & 0
    \end{array} \right]
\]
is the $6 \times 6$ nilpotent Weyr matrix of Weyr structure $x = (3,2,1)$.
\prskip

Another important reason for the success of the Weyr form in studying commuting problems (in addition to the critical property above) is that any matrix $K$ that centralizes (commutes with) a nilpotent Weyr matrix $W$ has an especially nice form, a blocked matrix version of the Toplitz matrices which centralize a basic nilpotent Jordan matrix $J$. This was first observed by Belitskii in 1983 (an account in English was given in \cite{B2}). We won't state the description precisely here, just illustrate with the above $W$ of structure $(3,2,1)$. Here a  centralizing matrix $K$ must look like
\[ \renewcommand{\arraystretch}{1.1}
            K \ = \  \left[ \begin{array}{r r r |r r |r }
    a & b  & d  &  g  & i  & l    \\
    0 & c  & e  &  h  & j  & m    \\
    0 & 0  & f  &  0  & k  & n    \\  \cline{1-6}
      &    &    &  a  & b  & g    \\
      &    &    &  0  & c  & h    \\  \cline{4-6}
  \multicolumn{5}{c |}{ }  & a
    \end{array} \right].
     \renewcommand{\arraystretch}{1.1}
\]
The centralizing matrices are then completely determined by their top row of blocks. We denote the $(i,j)$ block of a centralizing matrix $K$ by $K_{ij}$. Thus in the above displayed $K$, its top row blocks are
\[
   K_{11} \ = \ \left[\begin{array}{ccc}
                    a & b & d \\
                    0 & c & e \\
                    0 & 0 & f
                    \end{array}\right], \ \
  K_{12} \ = \ \left[\begin{array}{cc}
                    g & i \\
                    h & j \\
                    0 & k
                    \end{array}\right], \ \
 K_{13} \ = \ \left[\begin{array}{c}
                    l  \\
                    m \\
                    n
                    \end{array}\right].
\]
The forms of the top row blocks, outside of the homogeneous case, are not arbitrary but follow a certain staircase template. See \cite{HO}, \cite{O}, or \cite{ATLA} for more details. This approach is important for computing: given the Weyr matrix $W$, we construct the other pair $K,L$ of strictly upper triangular, commuting  matrices by specifying the centralizing form for their top row of blocks (thereby guaranteeing both commute with $W$), and then solving linear equations on their top row entries to ensure $K$ and $L$ commute with other. Such calculations done with the Jordan form would be extremely difficult (the centralizers don't take a nice block upper triangular form).
\prskip

\textbf{Summing Up}: \emph{In the reduction above to $W,K,L$ with $W$ in Weyr form, not only are $K,L$ strictly upper-triangular as ordinary matrices, they also have the nice centralizing form as block matrices relative to the Weyr structure of $W$. A win-win}!
\prskip

One of the on-screen, blow-by-blow, displays from our main routines is the output of so-called \textbf{\emph{leading edge dimensions}} \textbf{(\Small{LED})} resulting from the constructed commuting triple $W,K,L$. It is not an essential display but one of the most interesting, and possibly the best guide to how close we are to getting Eurekas for our particular setting of parameters. (Users who don't want this or other displays can easily suppress them.)  We won't go into the definition of the {\Small{LED}} here (the interested reader can look at \cite{O}, \cite{HO}, or \cite{ATLA}), but rather illustrate a couple of features by way of examples (taken from actual trials).
\thmskip

\noindent \textbf{Example 1.} \,Suppose we  are constructing commuting triples $W,K,L$ where $W$ has the Weyr structure $x = (8,8,8,4,1,1,1,1,1,1)$. Thus matrix size is $n = 8 + 8 + 8 + 4 + 1 + 1 + 1 + 1 + 1 + 1 = 34$. Here is an actual display  ($s$ is the seed used for the random number generator in {\Small{MATLAB}}):
\begin{align}
s \ &= \ 8146 \notag \\
 \dim \ &= \ 32 \notag \\
 \mbox{\Small{LED}} \ &= \ \begin{array}{|c|c|c|c|c|c|c|c|c|c|} \cline{1-10}
        8&8&8&4&1&1&1&1&1&1  \\ \cline{1-10}
        5&8&10&3&1&1&1&1&1&1   \\ \cline{1-10}
        \end{array} \notag
\end{align}
In the two-line display for the {\Small{LED}}, the first row gives the Weyr structure components. Beneath each is the corresponding leading edge dimension for that component. The output dim is $\dim F[W,K,L]$  provided  all steps were completed during the construction, but $\dim = 0$ indicates that certain linear equations had no solutions and the process had to stop. Now notice:
{\emph{
\begin{enumerate}
\item   $\dim$ is the sum of the {\Small{LED}}.
\item   If two successive Weyr structure components are the same, then the {\Small{LED}} for the second is at least as big as the first. (E.g. for the second and third components.)
\end{enumerate}
}}
These features hold in general.  In fact, the sum of the {\Small{LED}} at any given stage agrees with the dimension of the corresponding algebra of top-left corners (projecting here is an algebra homomorphism because the matrices are block upper-triangular.) So if the program is set up to display in real time (step-by-step), one can anticipate a Eureka if the process were to complete. This is illustrated in our second example.    \hfill $\square$
\prskip

\noindent \textbf{Example 2.}\, Consider the output
\begin{align}
s \ &= \ 8206 \notag \\
 \dim \ &= \ 0 \notag \\
  \mbox{\Small{LED}} \ &= \ \begin{array}{|c|c|c|c|c|c|c|c|} \cline{1-8}
        8&8&8&8&2&1&1&1  \\ \cline{1-8}
        6&9&9&0&0&0&0&0   \\ \cline{1-8}
        \end{array} \notag
\end{align}
Here the sum of the first three {\Small{LED}} is 24. And since the $4^{th}$ structure component is the same as the $3^{rd}$, if the $4^{th}$ step were to complete, we would have the sum of the {\Small{LED}} to this point at least $24 + 9 = 33 > 8 + 8 + 8 + 8 = 32$. This would give a Eureka for $F[\pi(W),\pi(K),\pi(L)]$ where $\pi$ is the projection map onto the top-left $4 \times 4$ corner of blocks.  Note that 0 beneath a structure component means the process did not get that far. \hfill $\square$
\prskip

In view of the above discussion, it is tempting to arrange that the $1^{st}$ {\Small{LED}} (= $\dim (\mathbb{Z}/p)[K_{11},L_{11}]$) be large. The biggest this can be is actually the $1^{st}$ structure component $n_1$ (by Gerstenhaber's theorem for commuting pairs --- note $W$ contributes nothing to the dimension of the first corner). A nontrivial theorem for commuting pairs of $k \times k$ matrices $A,B$ says that if $\dim F[A,B] = k$, then the only matrices that centralize both $A$ and $B$ are polynomials in $A,B$ (see \cite{LL}, \cite{NS}). From this it follows that in the homogeneous case (and in trials so far for the nonhomogeneous cases), if the $1^{st}$ {\Small{LED}} is the maximum possible (namely $n_1$), then $\dim F[W,K,L] = n$. So a Eureka can never occur this way. This illustrates an important principle: \emph{don't get greedy with the $1^{st}$ {\Small{LED}}}. On the other hand, if the $1^{st}$ {\Small{LED}} is small, there may be too much ground to make up to obtain a Eureka. It is a delicate balance.
\introskip

\section{The three main routines}
\thmskip

Here we outline the three main routines in our {\Small{MATLAB}} package. Various subroutines that the main routines call will be described in Sections 4, 5, 6, 7. There is a total of 38 routines in all (plus 3 sample trials)\footnote{\ The number of printed pages of code from our package runs to 63.}, but all most users will need to do is to run some of the main ones. Examples of this will be given in Section 8. Of course, if the user wants to modify some routines, that will require a better understanding of the program code. For those not very familiar with programming in {\Small{MATLAB}}, we recommend the text \cite{HH}.
\prskip

The three principal routines are \ \texttt{CommTriplesI,\, CommTriplesII,\, CommQuads}. We look at each in turn. (Note that \texttt{Comm} is an abbreviation for Commuting.)
\introskip

\noindent \texttt{CommTriplesI}.\, This has 15 input arguments and 5 output arguments:

\[
          \texttt{[EUREKAS,ES,W,K,L]  =  CommTriplesI(x,p,p0,p1,p2,pc,Min,Max,ICR,CS,TR,a,b,s1,s2)}
\]
\prskip

\noindent The most important inputs are the chosen Weyr structure \texttt{x} of the desired nilpotent Weyr matrix \texttt{W} (so $x$ is given as a partition of the chosen matrix size $n$, which is the sum of the partition entries, and so need not be separately inputted), the prime \texttt{p} for the mod $p$ arithmetic, and the seed range \texttt{[s1,s2]}. Thus the program runs through each seed $s$ in the integer interval $[s_1,s_2]$ and attempts to construct a commuting triple $A,B,C$ with $A = W$. If the process completes (more on that later), the routine calculates $\dim \,(\mathbb{Z}/p)[A,B,C]$. If this dimension exceeds the matrix size $n$, then the seed $s$ that produced it is recorded in the vector \texttt{ES} (of ``Eureka Seeds''). The first triple for which this happens is recorded in the outputs \texttt{W,K,L}. If no Eureka occurred during the run, all \texttt{W,K,L} are set to zero, and \texttt{ES} will be the empty array [ ]. On the other hand, the output \texttt{EUREKAS} returns the number of Eurekas that occurred (\texttt{EUREKAS = 0} indicates there were none.)\footnote{\ If \texttt{CommTriplesI} is being run continuously over several hours or days, it is useful to display the output \texttt{EUREKAS} periodically, so one can check for success from time to time.}
\prskip

The inputs \texttt{p0,\,p1,\,p2,\,pc} determine the sparsity of entries at various stages. Thus \texttt{p0,\,p1,\,p2} take values in the real interval $[0,1]$ and determine the probability of a particular entry being zero. The first of these gives sparsity for the first blocks $B_{11}, C_{11}$ of $B$ and $C$. After that sparsity is set as \texttt{p1} for the remaining set percentage \texttt{pc} of the blocks $B_{12}, C_{12}, \ldots $, then at \texttt{p2} for the rest of the blocks. The reason for this comes from experience with constructing commuting Eureka quadruples: \texttt{p1} is set lower so as to produce an early {\Small{LED}} that has jumped considerably above the corresponding Weyr structure component. Then \texttt{p2} is set much higher so as to improve the chances of the process completing (the linear equations that determine this are more likely to have a solution). See the earlier comments on the behaviour of the {\Small{LED}}.
\prskip

Unlike \texttt{CommTriplesII}, the routine \texttt{CommTriplesI} chooses for the first block $B_{11}$ of $B$ an  $n_1 \times n_1$ nilpotent Weyr matrix, where $x = (n_1,n_2,\ldots,n_r)$, subject to certain constraints. The advantage of this is that then the first block $C_{11}$ of $C$ has a known form because of the earlier centralizing description, and we can choose such a $C_{11}$ randomly (subject to the above constraints). It turns out that there is no loss of generality in making this choice of $B_{11}$ if $x$ has a homogeneous structure. But it is also known that this simplification is not always possible in the nonhomogeneous cases, whence $B_{11}$ and $C_{11}$ should be computed as random, strictly upper triangular, commuting $n_1 \times n_1$ matrices, as \texttt{CommTriplesII} does. This is no problem if $n_1$ is smallish, say $n_1 \le 10$. But it really slows the number of trials one can do in a given run once say $n_1 \ge 15$. It is a tradeoff between speed and generality. But since the analogue, \texttt{CommQuads}, of \texttt{CommTriplesI} works well in producing Eurekas for commuting quadruples of matrices, we tend to favour making the first corner of $B$ a Weyr matrix. The remaining blocks $B_{1j}, C_{1j}$ are then chosen from random solutions of certain linear equations.
\prskip

Inputs \texttt{Min} and \texttt{Max} give, respectively, the minimum and maximum allowable $1^{st}$ {\Small{LED}}. Setting \texttt{Min} $= (0.6)n_1$ and \texttt{Max} $= (0.9)n_1$ often work well, but in some situations the restrictions need to be stronger. Experimentation is advised.
\prskip

Input \texttt{ICR}\, (``Index Corner Restriction'') sets the proportion of $n_1 = x(1)$ (= $1^{st}$ corner size) to bound the nilpotent indices of $B_{11}$ and $C_{11}$. This imposes a bound (\texttt{ICR} $\times n_1 + r$) on the nilpotent indices of the final constructed $B,C$. See \cite{O} for the rationale for this. Setting \texttt{ICR} = 0.5 often works well.
\prskip

The input \texttt{CS}\, (``Change Step'') sets the number of trials run for the fixed choice of $1^{st}$ blocks $B_{11}$ and $C_{11}$, before changing to a new pair (but keeping the same $W$). And input \texttt{TR}\, ``Test Runs'' sets the number of trials in an optimizing process to determine the choice of the next random (but subject to the above constraints) pair $B_{11}, C_{11}$ so as to improve the completion rate for resulting triples $A,B,C$ (by increasing the rank of a certain coefficient matrix G in a linear system). Setting \texttt{TR} = 3000 often works well.
\prskip

Finally, the inputs \texttt{a,b} determine the interval of fit $[a,b]$ for a beta distribution for the output of dimensions.  One can take $b = (n + 1) + 0.5$ and $a$ to be 0.5 less than the theoretical minimum dimension (or estimate thereof). The accuracy of the estimate of $a$ is not critical. The fitted distribution is updated and displayed periodically, after every 1000 trials. Also displayed is the \texttt{estimate} so far of the probability of a Eureka occurring, \textbf{\emph{{\Small{IF}}}} one exists for our particular settings. On this basis, a user can decide whether to continue or try new settings. We stress, however, that the beta distribution is \textbf{\emph{{\Small{NOT}}}} used as a predictor of whether a Eureka exists here, but rather the number of trials needed to reasonably assure us otherwise. See \cite{O} for more details. In the case of commuting quadruples $A,B,C,D$ of matrices, our beta fit worked very well for the predicted estimates (say of getting $\dim (\mathbb{Z}/p)[A,B,C,D] = n+1$ or $n+2$, or even $n+3$). In large trials the estimates closely matched the actual number of these. So we see no reason why the predictions can't be relied on also for commuting triples.
\introskip

\noindent \texttt{CommTriplesII}.\, This has 11 input arguments, 9 in common with \texttt{CommTriplesI}, and the 6 output arguments, 5 in common with \texttt{CommTriplesI}:

\[
          \texttt{[EUREKAS,ES,W,K,L,HIST]  =  CommTriplesII(x,p,p1,p2,pc,B11,C11,a,b,s1,s2)}
\]
\prskip

\noindent The inputs that have been dropped are \texttt{p0, Min, Max, ICR, CS, TR}. The two new inputs are \,\texttt{B11,C11} and these are the first blocks of the constructed pair $B, C$. We give two subroutines for this (\texttt{CommPairI, CommPairII})\footnote{\ \texttt{CommPairII} is faster although not always guaranteed to produce strictly upper-triangular matrices.} , the user can choose one, that calculate  random, strictly upper-triangular, commuting $n_1 \times n_1$ matrices $B_{11}, C_{11}$. Users can, if they wish, achieve the earlier constraints involving Min and Max, and the size of the $1^{st}$ \texttt{LED}, by running the subroutines enough times to achieve these. Again, $x = (n_1,n_2,\ldots,n_r)$ is the Weyr structure of the first matrix $A$ of the commuting triples $A,B,C$ that is constructed. The new output \texttt{HIST} records the number of times $\dim = k$ occurred during the trials for $k = 1, \ldots ,n+5$. This enables a series of trials for different $B_{11},C_{11}$ to be run and their outputs to be collated (see \texttt{SampleTrial2}).
\introskip

\noindent \texttt{CommQuads}.\, This constructs for a given Weyr structure \texttt{x} and seed range $\texttt{[s1,s2]}$ commuting quadruples \texttt{W,K,L,M} of matrices with  \texttt{W} the nilpotent Weyr matrix of structure \texttt{x}:

\[
          \texttt{[EUREKAS,MaxExcess,MaxSeed]  =  CommQuads(x,p,p0,p1,p2,pc,a,b,s1,s2)}
\]
\prskip

\noindent The constructed quadruples are not displayed within the routine (they could be if the user wishes), although dimensions and {\Small{LED}} are displayed. The routine is designed principally to suggest promising parameter settings for snaring Eurekas in the case of commuting triples. The code for $\texttt{CommQuads}$ is less sophisticated than that for commuting triples. The inputs are a subset of those for $\texttt{CommTriplesI}$. Again, $\texttt{EUREKAS}$ is the number of Eurekas the trials produced. Output $\texttt{MaxExcess}$ is the maximum excess of dimension over matrix size that occurred in the production of Eurekas, while $\texttt{MaxSeed}$ is the first seed that gave MaxExcess. Another promising sign of a good choice of parameters is by how far some of the {\Small {LED}} have jumped above the corresponding Weyr structure component.
\introskip

\section{Three important supporting subroutines}
\prskip

The three most important subroutines are \texttt{WorkHorseI}, \;\texttt{WorkHorseII}, and \texttt{Dim}.
\thmskip

\noindent \texttt{WorkHorseI}.

\[
          \texttt{[dim,LED,E,S,W,K,L]  =  WorkHorseI(x,p,p1,p2,pc,K11,L11,G,C,PivInfo,s)}
\]
\prskip

\noindent This constructs a commuting triple \texttt{W,K,L} of matrices for a single inputted seed \texttt{s} and Weyr structure \texttt{x}. The first matrix is the nilpotent Weyr matrix of structure \texttt{x}. The inputs \texttt{p,\,p1,\,p2,\,pc} are the same as for \texttt{CommTriplesI}, while \texttt{K11,\,L11} are the $1^{st}$ corners of $K,\,L$. The flags \texttt{E,\,S} signal if a Eureka was produced, respectively whether the process completed, according to the values 1 and 0. Output \texttt{dim} is $\dim (\mathbb{Z}/p)[W,K,L]$ if the process completed, otherwise \texttt{dim = 0}. And \texttt{LED} gives the leading edge dimensions. To understand the inputs \texttt{G,\,C,\,PivInfo}, and only users who might want to modify the code need to read the details carefully, we need to go into how the top row blocks $K_{1k},\,L_{1k}$ are calculated for $k > 1$. Let $x = (n_1,n_2,\ldots,n_r)$.
\prskip

The $K_{1k},\,L_{1k}$ are calculated recursively. The construction of $X = K_{1k}$ and $Y = L_{1k}$ from the earlier top row blocks is determined by a linear system
\[
      Gv \ = \ w
\]
where
\[
      G \ = \ \left[\begin{array}{c}
                  M \\
                  D
              \end{array}\right],
\]
with $D$ a certain $2ab \times 2ab$ diagonal matrix for $a = n_1, \,b = n_k$,  and $M = [M_1,M_2]$  an augmentation of two $ab \times ab$ matrices
\begin{align}
   M_1 \ &= \ \ \ \ I_b \otimes K_{11} - K_{11}^t(1:b,1:b ) \otimes I_a \notag \\
   M_2 \ &= \ - \ (I_b \otimes L_{11} - L_{11}^t(1:b,1:b) \otimes I_a). \notag
\end{align}
The matrices $v, w$ are column vectors
\begin{align}
   v \ &= \ \left[\begin{array}{c} \mbox{vec}(Y) \\ \mbox{vec}(X)\end{array}\right], \notag \\
   w \ &= \  \left[\begin{array}{c}\mbox{vec}(Z) \\ 0 \end{array}\right],    \notag
\end{align}
where 0 is the $2ab \times 1$ zero matrix, and
\begin{align}
    Z \ &= \ K_{12}L_{1,k-1}(1:n_2,1:n_{k-1}) \,+ \; K_{13}L_{1,k-2}(1:n_3,1:n_{k-1}) \; + \; \cdots \ \notag \\ & \ \ \ + \ K_{1,k-1}L_{11}(1:n_{k-1},1:n_{k-1})  \ -\; L_{12}K_{1,k-1}(1:n_2,1:n_{k-1}) \notag \\
     & \ \ \ - \;L_{13}K_{1,k-2}(1:n_3,1:n_{k-1})\,- \; \cdots \ - \; L_{1,k-1}K_{11}(1:n_{k-1},1:n_{k-1}). \notag
\end{align}
Recall the vectorisation of an $m \times n$ matrix $X$, vec($X$), is the $mn \times 1$ column vector of the stacked columns of $X$.
\prskip

The system\; $Gv = w$\; comes from the requirement that the top left $k \times k$ corners of blocks of $K$ and $L$ commute (role of $M$) and that the new blocks conform to the shape of matrices centralizing $W$ (role of $D$). Notice that $G$ is $3ab \times 2ab$, \, $v$ is $2ab \times 1$, and $w$ is $3ab \times 1$. To solve the system we put $G$ in loose row-echelon form using elementary row operations, and record an invertible matrix $C$ such that $CG$ is the desired form (and so (CG)v = Cw is an equivalent system). Notice that each $G$ associated with the pair of blocks in a given position (there are $r - 1$ of these pairs) is completely determined by the first blocks $K_{11}, L_{11}$ and the structure $x$. Therefore, if we are performing a whole series of trials by running through a seed range $[s_1,s_2]$ (and thereby generating new solutions for the following blocks of $K$ and $L$), but involving the same first corners $K_{11}, L_{11}$, it makes sense to record all the $G$ and $C$ that are to be used. (That way, we avoid having to repeat the row operations each time.) This greatly increases efficiency of trialing (much faster).  When the individual $G$ are stacked in an array, this gives the input argument \texttt{G}. To make this possible, before stacking, we augment each $G$ with a zero matrix if necessary so that all the $G$ have the same number of columns, namely $n_1n_2$. Similarly, stacking the individual $C$, after augmenting with a zero matrix, gives the input \texttt{C}. From these arrays we can quickly recover each of the original, individual $G,C$ as required. Note that \texttt{G} is $3n_1(n - n_1) \times 2n_1n_2$ and \texttt{C} is $3n_1(n - n_1) \times 3n_1n_2$.
\prskip

This leaves only the input argument \texttt{PivInfo}. It is a $(\,3(r - 1) + 1) \times 2n_1n_2$ array that records the pivot positions (and more) of the various reduced matrices $CG$, thereby allowing rapid solving of the various linear equations as they arise.

\introskip

\noindent \texttt{WorkHorseII}.

\[
          \texttt{[dim,LED,E,S,W,K,L]  =  WorkHorseII(x,p,p1,p2,pc,K11,L11,s)}
\]
\prskip

\noindent This routine is similar to \texttt{WorkHorseI} but without the inputs \texttt{G,C,PivInfo}. They are calculated within the routine itself. The routine is designed to be used in conjunction with routines \texttt{CommPairI} or \texttt{CommPairII} to produce the $1^{st}$ corners \texttt{K11,L11}. Since this is much slower than the method used in \texttt{CommTriplesI}, to include  \texttt{G,C,PivInfo} as inputs as well would make \texttt{WorkHorseII} too slow, particularly over short runs.
\prskip

\texttt{WorkHorseII} can be used in its own right to produce a commuting triple of matrices. But its main purpose is to support the code of \texttt{CommTriplesII}. The subroutine is also used to produce the $1^{st}$ corners of $K,L,M$ in \texttt{CommQuad}. A similar technique could be used to construct $k$ commuting matrices $A_1,A_2,\ldots,A_k$, with the first matrix a nilpotent Weyr matrix, by  constructing the first corners of $A_2,A_3,...,A_k$ using the program for $k - 1$ commuting matrices. \texttt{WorkHorseIII} is the analogue of \texttt{WorkHorseII} for commuting quadruples.

\introskip

\noindent \texttt{Dim}.

\[
          \texttt{[dim,LED]  =  Dim(W,K,L,x,p)}
\]
\prskip

\noindent The subroutine \texttt{Dim} returns the dimension \texttt{dim} of the subalgebra generated by the three inputted commuting matrices \texttt{W,K,L}, where \texttt{W} is the nilpotent Weyr matrix of Weyr structure \texttt{x}. All calculations are done mod $p$ for the inputted prime \texttt{p}. Output \texttt{LED} gives the leading edge dimensions.
\introskip

Here is a sketch of how the subroutine works. Let $x = (n_1,n_2,\ldots,n_r)$. Suppose the nilpotent indices of $W,K,L$ are respectively $a,b,c$ (actually $a = r =$ \ length of the partition $x$). Then $(\mathbb{Z}/p)[W,K,L]$ is spanned as a vector space by the words
\[
       W^iK^jL^k \ \ \ \mbox{for} \ \ 0 \le i \le a - 1, \ 0 \le j \le b - 1, \ 0 \le k \le c - 1.
\]
Taking the vectorisation,\; $\mbox{vec}( W^iK^jL^k)$,\; of each of these words and making them the columns of a matrix $M$, we now have that dimension is the rank of $M$. There are a number of ways to find $\mbox{rank}(M)$. For instance, we could call the built-in routine in {\Small{MATLAB}} (it would need to be modified for mod $p$ arithmetic). However, this would be very slow because the matrix $M$ can be  huge: it is $n^2 \times abc$, where $n$ is the matrix size of $W,K,L$. Also, $a,b,c$ can be as large as $n$, so potentially we are looking at a $n^2 \times n^3$ matrix. In our search for Eurekas, we may have to consider matrices bigger than $100 \times 100$, potentially making $M$ bigger than $10,000 \times 1,000,000$. In practice we place restrictions on the nilpotent indices, usually not allowing indices bigger than 20. Still, $10,000 \times 8,000$ are big matrices if we are  running tens of thousands of trials.
\prskip

Our way around these problems is to take advantage of the special form the matrices $K,L$ have that comes from their centralizing of $W$. Firstly, since the $W^iK^jL^k$ centralize $W$, they are completely determined by their top row blocks. Thus we need only take the vectorisation of their top row blocks. We also order the set$ \{(i,j,k)\}$ of indices in a special way so that the corresponding parts of the matrix $M$ form an approximately lower triangular matrix. This allows us to calculate  $\mbox{rank}(M)$  recursively, using various elementary row operations to achieve a loose, row-echelon form of $M$ (whence $\mbox{rank}(M)$ = the number of nonzero rows of the reduced form). The earlier pivots are noted for subsequent clearing below but the earlier corner has zeros to its right in $M$ and is not touched again. Basically, we are recursively calculating the dimension of the algebra generated by the top left $m \times m$ corner of blocks of $W,K,L$ for $m = 1,2,\ldots, r$. We need these dimensions anyway for the {\Small{LED}}. Also we can check whether a Eureka has occurred at one of these corners. This is important because an early Eureka could be lost for the full matrices.
\prskip

In summary, our subroutine \texttt{Dim}  computes dimension very quickly and accurately.

\sectskip

\section{Subroutines connected with the Weyr form}
\introskip

We won't go into much detail here. The interested reader can look at the code and the accompanying documentation. The letter \texttt{W} in a routine name is shorthand for ``Weyr''.
\prskip

The routine \texttt{W = WMatrix(x)} returns the nilpotent Weyr matrix \texttt{W} of a given Weyr structure \texttt{x}, given as a partition $x = (n_1,n_2,\ldots,n_r)$ of the underlying matrix size $n$. It is independent of the underlying field. On the other hand \texttt{x = WStruc(A,p)} gives the Weyr structure \texttt{x} of a given nilpotent matrix \texttt{A} mod \texttt{p}. To complete the \emph{``three Weyr amigos''}, we have \texttt{[W,C] = WForm(A,p)}, which returns the Weyr form \texttt{W} of the inputted nilpotent matrix $A$, together with an invertible matrix \texttt{C} that achieves the similarity transformation ($C^{-1}AC = W$). The arithmetic is done modulo the prime \texttt{p}. But with all our routines, it would not be hard for a user to change the arithmetic to suit other situations (such as rational matrices), once safeguards are built in to avoid roundoff errors.
\prskip

We talked in Section 2 about how a matrix $K$ that centralizes a given nilpotent Weyr matrix $W$ of structure $x$ has a special blocked form, and is determined by its top row of blocks. In turn, these blocks have a certain staircase form. The routine \texttt{T = Template(x)} gives a template for these blocks, as an $n_1 \times n$ matrix of zeros and ones, where a 0 entry says that this entry must be zero, and an entry 1 says that the entry in this position can be arbitrary. This is often used in tandem with the routine \texttt{K = TopRowToFull(A,x)} to obtain the full $n \times n$ centralizing matrix $K$ that has the $n_1 \times n$ matrix \texttt{A} as its top row of blocks. Finally \texttt{K = RandCentMatrix(x,p,q,s)}  returns a random, subject to sparsity probability \texttt{q}, strictly upper-triangular matrix \texttt{K} that centralizes the nilpotent Weyr matrix $W$ of the inputted structure \texttt{x}. The answer depends on the inputted seed \texttt{s} that is to be used in the {\Small{MATLAB}} random number generator. The input \texttt{p} is the prime that determines the range $\{0,1,2,\ldots,p-1\}$ of allowable entries in \texttt{K}.
\sectskip

\section{Subroutines for solving linear equations}
\introskip

As is to be expected, if our routines involve solving linear equations, we should be using row operations. At the risk of reinventing the wheel, we have written our own versions to suit arithmetic mod $p$ and our particular inductive techniques in addressing the GP. Minimising the number of row operations (such as not making pivots 1 or swapping rows) and speed of execution are our priorities. The letter \texttt{L} in a routine name is shorthand for ``loose''.
\prskip

The routine \texttt{LFormI} takes the ``form''  \texttt{[L,C,P,r] = LFormI(A,TI,p).} Here \texttt{A} is the matrix to be reduced and \texttt{p} is the prime for mod arithmetic. The output \texttt{L} is  row equivalent to \texttt{A} and has a loose row echelon form: no two pivots (first nonzero entries of nonzero rows) lie in the same column, pivots need not be 1, and zero rows can be interspersed with nonzero rows. Note that our form need not have a staircase (echelon) shape. Output \texttt{r} is the rank of \texttt{A} (the number of nonzero rows of \texttt{L}), and \texttt{P} is a column vector that records the pivot positions: \texttt{P(j) = i} if row $i$ contains a pivot in column $j$, otherwise \texttt{P(j) = 0}. Input \texttt{TI} (``To Invert'') takes the values 1 and 0 according to whether we wish to produce an invertible matrix that gives (under left multiplication) the form \texttt{L} or not (for instance, we don't need the invertible matrix if all we are interested in is the rank of $A$). If so desired, output \texttt{C} gives such an invertible matrix, otherwise $C$ is the identity matrix.
\prskip

Routine \texttt{[L,C,r] = LFormII(A,TI,p)} is similar to \texttt{LFormI} except now output \texttt{P} has been dropped, and entries above and below a pivot are 0 (whence $C$ may no longer be lower triangular).
\prskip

Routine \texttt{[L,C,Q,r] = LFormIII(A,B,TI,P,r1,p)} is more interesting. Here the input \texttt{A} is an $m \times n$ matrix already in the form of \texttt{L} in \texttt{LFormI}, \texttt{P} is the column vector of the pivot positions, and \texttt{r1} is the rank of \texttt{A}. The input \texttt{B} is a $t \times n$ matrix. Then output \texttt{L} is a matrix in the loose form of \texttt{LFormI} of the augmented matrix
\[
       \left[\begin{array}{c}
                A \\
                B
                \end{array}\right].
\]
Output \texttt{Q} gives the new pivot positions, \texttt{r} is the new rank, and \texttt{C} is an invertible matrix that realises the reduction to \texttt{L}.
\prskip

The routine \texttt{Inv} uses \texttt{LFormII} to invert a matrix or record that it is not invertible. The techniques in our routines for solving a linear system $Ax = b$ using some of the above forms are pretty standard. \texttt{SolveI} uses \texttt{LFormI} and works through the equations in the column-order of the pivots, starting with far rightmost pivot and working backwards. On the other hand, \texttt{SolveII} uses \texttt{LFormII} and reads off the solutions directly.
\sectskip

\section{The beta fit routine}

\introskip

The final subroutine we have chosen to highlight is \texttt{BetaFit}. As described in Section 3 it is a useful guide to the number of trials needed to flush out a Eureka \textbf{\emph{if}} such exist for our particular setting of parameters. A beta distribution takes the form
\[
     p(x) \ = \ Cx^d(1-x)^e, \ \ x \in [0,1]
\]
for some constants $C,d,e$ and over the real interval $[0,1]$ --- so it is a continuous distribution. (Formal descriptions state this differently in terms of ``shape parameters'' $\alpha,\beta$.)  On a general (finite) interval $[a,b]$, we need to replace $x$ in the right hand side of $p(x)$ by $(x-a)/(b-a)$ and replace $C$ by $C/(b-a)$ \,(but retain same $d,e$).
\prskip

From a series of trials for  commuting triples $W,K,L$ of $n \times n$ matrices, we view the histogram of the computed dimensions as a crude approximation to the true density function. (The histogram rectangle, say at $\dim = 10$, should be viewed as sitting on $[9.5,\, 10.5]$.) We then fit a beta distribution on $[a,b]$  to the histogram and from that estimate the probability of getting $\dim = n+1$ \,if such a Eureka exists.
\prskip

Figure 1 is the beta distribution (red)
\[
        p(x) \ = \ 414.66[((x-7.5)/14))^{8.5642}]\,[(1-((x-7.5)/14))^{4.0426}]
\]
for the histogram (blue) of dimensions from commuting triples of $20 \times 20$ matrices $W,K,L$ with $W$ the nilpotent Weyr matrix of structure $x = (7,5,5,2,1)$, and working mod the prime $p = 101$. (This was suggested by commuting quadruples which gave the promising output of dimensions 21,22, and even 23 with probabilities  0.0367, \,0.0039, \,$3.6297 \times (10)^{-5}$ respectively, which closely matched the actual number of occurrences.) The program used to produce the underlying data needed for Figure 1 is similar to that in \texttt{SampleTrial1} given later.  The green base is a reminder of the interval of fit. And the green $\ast$ are simply a visual guide as to where the base (integer) dimension is centred. The fit here for commuting triples looks good over the distribution of our gathered raw data, and suggests no reason why it can't be relied on for a Eureka probability prediction if one exists here.
\prskip

\prskip

On the assumption that  $p(x)$ gives even a halfway  reasonable fit to the upper tail of the underlying distribution, $\dim = 21 = n + 1$ should be popping up about 1 in 500 trials (and certainly at least every 10,000 trials). The fact that it didn't occur at all in some 2 million trials makes it a reasonable inference that the GP does not fail here. So we should look elsewhere.
\prskip

Here is the form of the routine \texttt{BetaFit}:
\[
\texttt{[alpha,beta,C,M,V,Estimate] = BetaFit(n,a,b,Hist)}
\]
Most of the input and output arguments are self-explanatory. Input \texttt{Hist} is the histogram of the dimensions recorded from a (large) series of trials. It takes the form of a $1 \times n$ row vector whose $i^{th}$ component is the number of times $\dim = i$ occurred. ($n$ as always is the underlying matrix size.) Output \texttt{Estimate} gives the probability of a Eureka occurring if one exists. The calculations used to fit the beta probability density function involve only the mean and variance of the data, given as outputs \texttt{M,V}, not higher moments. But it works fine because we are not after a precise probability estimate (that would be very difficult). The interested reader can look at the actual code for more detail.
\prskip

This concludes our description of some of the routines in our package. There are a total of 38 routines in all, although some of them are at quite a low level of programming. To assist the reader in keeping track of the overall structure, we include the following chart of dependency relations between certain highlighted routines. We use the notation
\[
{\SMALL
\xymatrix{ A \\
B \ar[u] }}
\]
to indicate program A calls program B.
\introskip

\noindent {\Small \textbf{DEPENDENCY}}
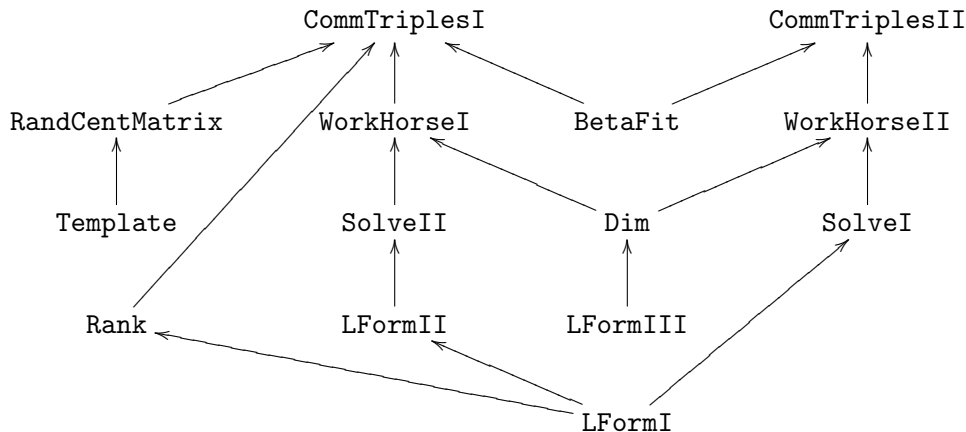
\introskip

\[
\xymatrix{  & \texttt{CommTriplesI} & & \texttt{CommTriplesII} & \\
 \texttt{RandCentMatrix} \ar@{>}[ru] & \texttt{WorkHorseI} \ar[u] &  \texttt{BetaFit}  \ar@{>}[ru] \ar@{>}[lu]  & \texttt{WorkHorseII}  \ar[u] & \\
\texttt{Template} \ar[u] & \texttt{SolveII} \ar[u] & \texttt{Dim} \ar@{>}[lu] \ar@{>}[ru] & \texttt{SolveI} \ar[u] & & \\
 \texttt{Rank} \ar@{>}[ruuu] & \texttt{LFormII} \ar[u] &  \texttt{LFormIII} \ar[u] &  &  & \\
& & \texttt{LFormI}  \ar@{>}[llu] \ar@{>}[lu] \ar@{>}[ruu] & & }
\]
\sectskip

\section{Some examples of running the main routines}
\introskip

Many Mathematics Departments have a general licence for {\Small{MATLAB}} that can be shared with the faculty. So for those not very familiar with {\Small{MATLAB}}, our advice is to enquire whether your department has such a licence.\footnote{\ Ask a friendly technician to install {\Small{MATLAB}} on your office computer and to give a brief tutorial on how to run programs.} Our package of programs for the Gerstenhaber Problem should then be downloaded into a folder within {\Small{MATLAB}}.
\prskip

Below are three sample trials involving actual {\Small{MATLAB}} code (they are also included in our package). Rather than entering the code directly into {\Small{MATLAB}}'s command window, it is better to create a separate m-file \texttt{SampleTrial.m} and execute that
by typing in the command window
\[
       \texttt{tic, SampleTrial, toc}
\]
and hitting the enter button. An advantage of this is that it is easy to make small changes (changing parameter settings, for example) in future trials. Note that the \texttt{tic, toc} will return the running time for the trial (these are optional).
\prskip

As a guide to the running times of each sample trial, we give the times from a basic 4-core-processor laptop of around 2.5 GHz speed (and using a 64 bit version of {\Small MATLAB}). Our first trial \texttt{SampleTrial1} illustrates \texttt{CommTriplesI} and takes around 1 hour.\footnote{\ The authors would regard this as very short trial.}
\introskip

\noindent
\texttt{\% SampleTrial1  \ \ \       Tests CommTriplesI on matrix size n = 29 \\
\% \hspace*{30mm}                                 with structure x = [9,7,7,4,1,1].\\
n = 29 \\
x = [9,7,7,4,1,1]  \\
p = 97  \ \ \ \ \ \ \ \% chosen prime;\\
p0 = 0.5; \\
s1 = 1  \\
s2 = 25000  \ \ \     \% so 25,000 trials; \\
a = 12.5; \\
b = 30.5; \\
ICR = 0.6 \\
CS = 1000;  \ \ \                 \% K11,L11 change after 1000 trials; \\
TR = 3000;  \ \ \                \% 3000 experiments for optimal K11,L11;  \\
p1 = 0.88  \\
p2 = 0.93  \\
pc = 50;     \ \  \  \  \           \% 2nd and 3rd blocks have p1, last three p2; \\
Min = 4 \\
Max = 7 \\
$[$EUREKAS,ES,W,K,L$]$ = CommTriplesI(x,p,p0,p1,p2,pc,Min,Max,ICR,CS,TR,a,b,s1,s2),
}
\thmskip

Our second trial illustrates \texttt{CommTriplesII} and takes under 1 and 1/2 hours.
\thmskip

\noindent
\texttt{\%  SampleTrial2 \ \ \      Tests CommTriplesII on matrix size n = 18 \\
\% \hspace*{30mm}                  with structure x = $[$7,4,3,3,1$]$. Makes 20 \\
\% \hspace*{30mm}                  changes in 1st corner in a series of runs \\
\% \hspace*{30mm}                   of 2000 trials each. Then fits beta \\
\% \hspace*{30mm}                   distribution (and estimates pr(dim = n+1))  \\
\% \hspace*{30mm}                   over all 40,000 trials. \\
n = 18 \\
x = $[$7,4,3,3,1$]$ \\
p = 37 \\
q = 0.82 \\
p1 = 0.75 \\
p2 = 0.90 \\
pc = 60 \\
a = 5.5 \\
b = 19.5  \\
EUREKAS = 0; \\
FirstES = 0; \\
Hist = zeros(1,n+5); \\
for t = 0:19 \\
\hspace*{6mm}     if t > 0  \\
\hspace*{12mm}           if EUREKAS == 0 \\
\hspace*{18mm}                  disp('No Eurekas so far.'),  \\
\hspace*{12mm}           else \\
\hspace*{18mm}                  disp('EUREKA WAS FOUND AT SEED'), FirstES \\
\hspace*{12mm}           end \\
\hspace*{12mm}           pause(5) \\
\hspace*{6mm}     end \\
\hspace*{6mm}    $[$B11,C11$]$ = CommPairII(x(1),1,p,q,t), \% could use CommPairI here instead; \\
\hspace*{6mm}     s1 = 1+2000*t; \\
\hspace*{6mm}     s2 = 2000 + 2000*t; \\
\hspace*{6mm}     $[$Eurekas,ES,W,K,L,HIST$]$ = CommTriplesII(x,p,p1,p2,pc,B11,C11,a,b,s1,s2); \\
\hspace*{6mm}     Hist = Hist + HIST; \\
\hspace*{6mm}     if ( (EUREKAS == 0) \& (Eurekas \textasciitilde 0) ) == 1   \\
\hspace*{12mm}            First ES = ES(1); \\
\hspace*{6mm}     end \\
\hspace*{6mm}     EUREKAS = EUREKAS + Eurekas; \\
end \\
EurekasOverAll = EUREKAS, \\
Hist, \\
display('Beta distribution over all our choices of B11,C11'), \\
Hist = $[$Hist(1:n),0,0,0$]$; \\
$[$alpha,beta,CC,MM,VV,Estimate$]$ = BetaFit(n,a,b,Hist); \\
FinalEstimate = Estimate, \\
TotalCompletedCases = sum(Hist), \\
}
\introskip

Our final trial illustrates \texttt{CommQuads}. Its code is a couple of pages long, so we refer the reader to the code and documentation given within our package of routines. The running time is around 2 and 1/4 hours. However, a user need not run the full trial to see what is going on because there are periodic updates on the results so far. The best way of getting an overview of the output is to watch for the changes in the displayed beta fit graph, at about 3 minute intervals after the completion of 1000 trials for a particular choice of settings. From this one can quickly see the number of Eurekas, maximum excess of dimension over matrix size, etc. Here is the introductory description of the trial.
\introskip

\noindent
\texttt{\% SampleTrial3 \ \ \   For matrix size n = 15, the trial illustrates  \\
\% \hspace*{30mm}               how to choose promising parameter settings for  \\
\% \hspace*{30mm}               CommTriplesI and CommTriplesII to snag Eurekas,  \\
\% \hspace*{30mm}               based on the routine CommQuads and promising \\
\% \hspace*{30mm}               settings for commuting quads of n x n  matrices \\
\% \hspace*{30mm}               that actually produce lots of Eurekas.
}

\introskip

The references below include those mentioned explicitly in our text along with a number of others that are related to the GP.

\introskip

\noindent {\Small{\textbf{Acknowledgement.} The authors wish to thank Ian Coope, Jean-Pierre Schoch, Pace Nielsen, Robert Carlson, and Graham Wood for their helpful advice.}}
\thmskip

\introskip

\noindent {\Small {\textbf{AFFILIATIONS}}}
\prskip

\noindent (John Holbrook) \ \ \ University of Guelph, Guelph, Ontario, Canada.\\
\noindent (Kevin O'Meara) \ \ University of Canterbury, Christchurch, New Zealand.
\introskip

\noindent {\Small {\textbf{E-MAIL}}} \\
\noindent \texttt{jholbroo@uoguelph.ca} \ \ \ (John Holbrook) \\
\noindent \texttt{staf198@uclive.ac.nz} \ \ \ (Kevin O'Meara)

\end{document}